\documentclass[oneside]{amsart}
\usepackage{amsmath}
\usepackage{amsfonts}
\usepackage{amssymb}
\usepackage{latexsym}
\usepackage{amsthm}
\usepackage{setspace}
\usepackage{hyperref}
\usepackage{geometry}
\usepackage[usenames,dvipsnames]{color}

\everymath{\displaystyle}
\usepackage{url}

\newcommand{\Oc}{\mathcal{O}}

\newcommand{\fm}[4]{\begin{pmatrix} #1 & #2 \\ #3 & #4 \end{pmatrix}}           

\newcommand{\Z}{\mathbb Z}
\newcommand{\Q}{\mathbb Q}
\newcommand{\A}{\mathbb A}

\newcommand{\F}{\mathbb F}

\newcommand{\C}{\mathbb C}
\newcommand{\p}{\mathfrak p}

\newcommand{\Ac}{\mathcal{A}}
\newcommand{\Gc}{\mathcal{G}}
\DeclareMathOperator{\GL}{GL} \DeclareMathOperator{\SL}{SL}  
\DeclareMathOperator{\Gal}{Gal} \DeclareMathOperator{\Ind}{Ind}

\newcommand{\Nrm}{\operatorname{\textbf{N}}}
\newcommand{\St}{\operatorname{\textrm{St}}}
\newcommand{\BC}{\operatorname{\textrm{BC}}}
\newcommand{\Res}{\operatorname{\textrm{Res}}}
\newcommand{\Art}{\operatorname{\textrm{Art}}}
\newcommand{\PS}{\operatorname{\textrm{PS}}}
\newcommand{\Ss}{\operatorname{\textrm{S}}}
\newcommand{\Dim}{\operatorname{\textrm{dim}}}
\newcommand{\Hom}{\operatorname{\textrm{Hom}}}
\newcommand{\tr}{\operatorname{\textrm{tr}}}

\newtheorem{thm}{Theorem}

\newtheorem{dfn}[thm]{Definition}
\newtheorem{prop}[thm]{Proposition}

\theoremstyle{remark}
\newtheorem{rem}[thm]{Remark}

\begin{document}

\title{Dimension Formulae for Spaces of Lifted Bianchi Modular Forms}

\author{Mehmet Haluk \c{S}eng\"un}
\address{Mathematics Institute, Zeeman Building, University of Warwick, Coventry, CV4 7AL, UK}
\email{M.H.Sengun@warwick.ac.uk}

\author{Panagiotis Tsaknias}
\address{Facult� des Sciences, de la Technologie et de la Communication, Universit\'e du Luxembourg}
\email{panagiotis.tsaknias@uni.lu}
\email{p.tsaknias@gmail.com}


\begin{abstract}
In this note, we work out the dimension of the subspace of base-change forms and their twists inside 
the space of Bianchi modular forms with fixed (Galois stable) level and weight. The formulae obtained here 
have been used for the numerical experiments of \cite{bsv}. 
Extended versions of the formulae obtained in this note, treating CM-forms and spaces with non-trivial nebentypus, will be employed in the upcoming paper \cite{even}.
\end{abstract}

\maketitle


\section{Setting}

Let $K=\Q(\sqrt{-D})$ be an imaginary quadratic field. Let $f\in S_k(\Gamma_0(N))$ be classical newform. To such an object one attaches an irreducible automorphic representation $\pi:=\pi_f$ of $\GL_2(\A_\Q)$. One then has the so-called \emph{Base Change} operator:
$$\BC_\Q^K:\Ac_n(\Q)\to \Ac_n(K)$$
We would like to understand how the base change operator affects the level of an automorphic representation. We will in fact restrict our attention to the $2$-dimensional case (which will in turn force us to consider the $1$-dimensional case as well).

Let $\pi\in\Ac_n(\Q)$. One can write it as a restricted tensor product of local ones $\pi_p\in\Ac_n(\Q_p)$, where  $p$ runs through all finite primes and infinity.  It is well-known that Base-Change is compatible with taking local components, i.e for every rational prime $p$ and  every prime $\mathfrak{p}$ of $K$ above it we have: 
$$\BC_\Q^K(\pi)_{\mathfrak{\p}}=\BC_{\Q_p}^{K_\mathfrak{p}}(\pi_p)$$
In light of this it is enough to examine what happens at each prime individually. This allows us to freely exploit the local Langlands correspondence throughout this note. In fact the base change operator is the one that corresponds on the restriction operation (locally and globally)
$$\Res_{W_K}^{W_\Q}: \Gc_n(\Q) \to \Gc_n(K),$$
$$\Res_{W_{K_\mathfrak{p}}}^{W_{\Q_p}}: \Gc_n(\Q_p) \to \Gc_n(K_\mathfrak{p}),$$
on $n$-dimensional Weil-Deligne representations. We will denote by $\Art_F:W_F\to F^\times$ the local Artin map of local class field theory. Then the local Langlands correspondence in the $1$-dimensional case is composition with $\Art_F^{-1}$. 

Going back to the automorphic side, if $p$ is finite then one has the following classification of the types of the $\pi_p$ that can occur:
\begin{itemize}
\item Principal series $\PS(\chi_1, \chi_2)$, where $\chi_i\in\Ac_1(\Q_p)$ and $\chi_1/\chi_2\neq ||.||^{\pm1}$. The corresponding Weil-Deligne representation $WD(\pi_p)=(\rho_p, N_p)$ is in this case split (decomposable), i.e. $\rho_{p}\sim \tilde\chi_1\oplus\tilde\chi_2$ and $N_p=0$, where $\tilde\chi_i$ is the corresponding character of the Weil Group, i.e. $\tilde\chi=\chi\circ\Art_p$ . We will study this type in Section \ref{sec:split}.
\item Special $\Ss(\chi||.||,\chi)$, where $\chi\in\Ac_1(\Q_p)$. The corresponding local Galois representation is in this case reducible but not decomposable: $\rho_{\pi_p}\sim\tilde\chi||.||\oplus\tilde\chi$ and $N=\fm{0}{1}{0}{0}$. This will be the subject of Section \ref{sec:special}.
\item Supercuspidal, which is all the rest. We will describe this type in more detail in section Section \ref{sec:cusp}.
\end{itemize}

Let now $p$ be a rational prime and $\mathfrak{p}$ a prime in $K$ above it. We have the following:
\begin{itemize}
\item $D_p=D_\mathfrak{p}$ iff $p$ splits, i.e. $p=\mathfrak{p}\mathfrak{p'}$.
\item $I_p=I_\mathfrak{p}$ iff $p$ is unramified, i.e $p\neq \mathfrak{p}^2$.
\end{itemize}
This trivial observation tells us that $\alpha(\rho|_{I_\mathfrak{p}})=\alpha(\rho|_{I_p})$ for all the unramified primes $p$ in $K$ and therefore the conductor stays the same at these primes. Lets examine what happens for the ramified primes in $K$ whose (finite) set we will denote by $R_K$.

\begin{rem} Since we have established the finiteness of the set of primes where any change in the conductor can possibly occur one can wlog consider a suitable $\ell$-adic representation of $G_\Q$ instead of the Weil-Deligne representation of $W_\Q$: This follows from the Deligne-Grothendieck classification of $\ell$-adic representations if we pick $\ell$ big enough so that $\ell$ not in $R_K$ and $\pi$ is unramified at it. Same thing goes locally.\end{rem}

In what follows, we will assume for simplicity that $2$ is unramified in $K$. We treat first the case $n=1$, and then in the following three sections we treat the case $n=2$, one type at a time. We then finish with some dimension formulae for spaces of modular forms of a given type. In the future, we will also consider the Automorphic Induction operator and spaces with non-trivial Nebentypus.


\section{The case $n=1$}

We will work out this case not just for quadratic extensions of $\Q_p$ but for any quadratic extension $E/F$ of local fields over $\Q_p$. Let $\chi\in\Ac_1(F)$. Let $\mathfrak{p}_F$  and $\p_E$ be the primes in $\Oc_F$ and $\Oc_E$ respectively such that $\p_E|\p_F|p$. We then have an explicit description of the Base Change operator thanks to local class field theory:
$$\BC_{F}^{E}(\chi)=\chi\circ\Nrm_{E/F}.$$

Recall that the conductor of a character $\eta$ of $F^*$ is the smallest integer such that $\eta|_{U_F^{i}}$ is trivial, with $U_F^{i}$ defined as:
$$U_F^{i}=\begin{cases}
\Oc_F^\times&i=0\\
1+\p_F^i&i>0\\
\end{cases}$$

Clearly if $\chi$ is unramified, then $\BC_F^E(\chi)$ is unramified too. Assume that $\chi$ is tamely ramified (and therefore $p>2$). It then factors through the quotient $\Oc_F^*/U^1_F\cong k_F^*$. By \cite[Lemma 18.1]{BushnellHenniart2006} we have that $\BC_{\Q_p}^{K_\mathfrak{p}}(\chi_p)$ can be at most tame. It is unramified if and only if $\chi_p|_{\Z_p^*}$ has order $2$.

The following Proposition explains what happens when the conductor of $\chi_p$ is at least $2$ in the case where $K_\mathfrak{p}/\Q_p$ is tamely ramified (which is the only possible case since $p>2$):
\begin{prop}{\cite[Proposition 18.1]{BushnellHenniart2006}}
Assume that $K_\mathfrak{p}/\Q_p$ is tamely ramified. Then:
$$\alpha(\BC_{\Q_p}^{K_\mathfrak{p}}(\chi_p))=e_{K_\mathfrak{p}/\Q_p}(\alpha(\chi_p)-1)+1\qquad \forall m\geq2$$
\end{prop}

\begin{rem} Bushnell and Henniart use the term \emph{level}, which is equal to the conductor minus $1$.\end{rem}

Putting everything together:
\begin{prop}
$$\alpha(\BC_{\Q_p}^{K_\mathfrak{p}}(\chi_p))=\begin{cases}
0&\alpha(\chi_p)=0\\
0&\alpha(\chi_p)=1\textrm{ and }\chi_p^2|_{Z_p^*}=1\\
1&\alpha(\chi_p)=1\textrm{ and }\chi_p^2|_{Z_p^*}\neq1\\
e(K_\mathfrak{p}/\Q_p)(m-1)+1&\alpha(\chi_p)=m\geq2\\
\end{cases}$$
\end{prop}

\begin{rem} Notice that this proposition covers the conductor at all primes of $K$, not just the ramified ones. As we expected, the conductor does not change at the unramified ones though.\end{rem}

\section{The case $n=2$}

\subsection{The split case}\label{sec:split}

In this case $\pi_p$ is a principal series $I(\chi_1,\chi_2)$ parametrized by two Hecke characters $\chi_1,\chi_2\in\Ac_1(\Q_p)$. Obviously one has (by looking at the Galois side for example) that
$$\BC_{\Q_p}^{K_\mathfrak{p}}(I(\chi_1,\chi_2))=I(\BC_{\Q_p}^{K_\mathfrak{p}}(\chi_1), \BC_{\Q_p}^{K_\mathfrak{p}}(\chi_2)).$$
Combining this with the conductor formula
$$\alpha(I(\chi_1,\chi_2))=\alpha(\chi_1)+\alpha(\chi_2)$$
for principal series one then reduces the problem to the $1$-dimensional case that we just treated in the previous Section.

\subsection{The Special case}\label{sec:special}

In this case $\pi_p$ is a (twisted) Steinberg representation $S(\chi||.||, \chi)=\chi\otimes\St$ parametrized by one Hecke character $\chi\in\Ac_1(\Q_p)$. We have the following conductor formula:
$$\alpha(\chi||.||, \chi)=\begin{cases}
1&\alpha(\chi)=0\\
2\alpha(\chi)&\alpha(\chi)>0.\\
\end{cases}$$
In order to understand base change in this case we will exploit the Galois side. The base change corresponds to restriction there and therefore it is easy to see that:
$$\BC_{\Q_p}^{K_\mathfrak{p}}(S(\chi||.||, \chi))=S(\BC_{\Q_p}^{K_\mathfrak{p}}(\chi)||.||, \BC_{\Q_p}^{K_\mathfrak{p}}(\chi))$$
Since the same dimension formula applies for the RHS, we have again reduced the problem to one in the $1$-dimensional case.

\subsection{The Supercuspidal case}\label{sec:cusp}

In this case $\pi_p$ is none of the previous types. Since we have the ongoing assumption that $p>2$, all these local automorphic representations are parametrized by pairs $(E,\theta)$ where $E/\Q_p$ is quadratic and $\theta\in\Ac_1(E)\backslash \BC_{\Q_p}^E(\Ac_1(\Q))$. We also get that the corresponding Galois representation looks like
$$\rho_{\pi_p}\sim\Ind_{G_E}^{G_{p}}(\theta\circ \Art_p).$$
We will write $\tilde\theta$ for $\theta\circ\Art_p$. Base change corresponds to restriction on the Galois side so we want to examine the conductor of $\Res_{G_\mathfrak{p}}^{G_p}\Ind_{G_E}^{G_{p}}\theta$.
Using \cite[Proposition 22]{Serre1977} we get:
$$\Res_{G_\mathfrak{p}}^{G_p}\Ind_{G_E}^{G_{p}}\tilde\theta=\begin{cases}
\tilde\theta\oplus\tilde\theta^\sigma&G_\mathfrak{p}=G_E,\\
\Ind_{G_E\cap G_\mathfrak{p}}^{G_\mathfrak{p}}\tilde\theta|_{G_E\cap G_\mathfrak{p}}&G_\mathfrak{p}\neq G_E\\
\end{cases}.$$
Here $\sigma$ is a generator of $\Gal(K_\mathfrak{p}/\Q_p)$ and $\tilde\theta^\sigma(x)=\tilde\theta(g^{-1}xg)$, for all $x\in G_\mathfrak{p}$ and $g\in G_p$ a lift of $\sigma$. Notice that since $K_\mathfrak{p}/\Q_p$ is ramified we have that $I_p\not\subseteq G_\mathfrak{p}$ so we can pick $g$ to be in $I_p$.

Let's treat the first case. The conductor of the representation in the Galois side is equal to $\alpha(\tilde\theta) + \alpha(\tilde\theta^\sigma)$. We also have
that $\alpha(\tilde\theta) = \alpha(\tilde\theta^\sigma)$ so we get that:
$$\alpha(\BC_{\Q_p}^{K_\mathfrak{p}}(\pi_p))=2\alpha(\theta).$$

For the second case we employ \cite[Theorem 8.2]{Henniart1979} which gives the following formula:
$$\alpha(\Ind_{G_{EK_\mathfrak{p}}}^{G_\mathfrak{p}}\tilde\theta|_{G_{EK_\mathfrak{p}}}) = f_{EK_\mathfrak{p}/K_\mathfrak{p}}\left(d_{EK_\mathfrak{p}/K_\mathfrak{p}} + \alpha(\tilde\theta|_{G_{EK_\mathfrak{p}}})\right)$$
Translating it into a statement on the automorphic side gives:
$$\alpha(\BC_{\Q_p}^{K_\mathfrak{p}}(\pi_p)) = f_{EK_\mathfrak{p}/K_\mathfrak{p}}\left(d_{EK_\mathfrak{p}/K_\mathfrak{p}} + \alpha(\BC_{K_\mathfrak{p}}^{EK_\mathfrak{p}}(\theta))\right)$$
The conductor on the RHS can then be computed using the formula obtained for the $1$-dimensional case.


\section{Dimension formulae for spaces of fixed type}
Let $N$ be a positive integer. We will denote the automorphic forms of level $N$ by $\mathcal{A}^N_2(\Q)$, i.e. $\mathcal{A}_2^N(\Q):=\mathcal{A}_2(\Q)^{U_1(N)}$. In the same fashion, $\mathcal{A}_2^{\mathfrak{n}}(K):=\mathcal{A}_2(K)^{U_1(\mathfrak{n})}$ will denote the automorphic forms over $K$ of level $\mathfrak{n}$, where $\mathfrak{n}$ is an ideal of $\mathcal{O}_K$.

Assume for simplicity that $K$ be an imaginary quadratic field with {\bf odd} discriminant $d=disc(K /\Q)$ and that $N$ is a {\bf square-free} positive integer that is coprime to $d$. Let $\Pi \in \mathcal{A}^{(N)}_{2,bc}(K)$ and $\pi \in \mathcal{A}_2(\Q)$ such that $\Pi$ is a twist of $BC_\Q^K(\pi)$. 

\begin{dfn} For $k \ge 0$, let $\mathcal{B}^N_2(\Q,k)$ denote the set of all cuspidal automorphic representations $\pi$ of $\GL_2(\mathbb{A}_\Q)$ 
such that
\begin{itemize}
\item $\pi_\infty$ is the holomorphic discrete series representation of weight $k+2$, 
\item for every prime $p | N$, $\pi_p$ is of Steinberg type,
\item for every prime $\ell | d$, there are three possibilities for $\pi_\ell$:
\begin{enumerate}
\item unramified principal series,
\item $PS(\alpha, \omega_{K,\ell} \beta)$ with $\alpha, \beta$ unramified characters of $\Q_\ell^*$ and $\omega_K$ is the quadratic character associated to $K$,
\item $AI_{K_\ell / \Q_\ell}(\theta_\ell) \otimes \gamma$ with an unramified chracter $\gamma$ of $\Q_\ell^*$.
\end{enumerate}
\end{itemize}
\end{dfn}

Assume that $\Pi$ is of weight $k$, that is, $\Pi_\infty$ is the principal representation $\rho^k_\infty=PS(\chi_\infty^{k+1} ,\chi_\infty^{-k-1})$ of $\GL_2(\C)$ with 
$\chi_\infty (z) = z/|z|$. Then one can find $\pi \in \mathcal{B}_2^N(\Q, k)$ such that $\Pi \simeq BC(\pi) \otimes \chi$ for some finite order idele class character $\chi$ 
of $K$. Conversely, for any non-CM $\pi$ in $\mathcal{B}_2^N(\Q,k)$, there exists a finite order idele class character $\chi$ 
of $K$ such that $\Pi \simeq BC(\pi) \otimes \chi$. Moreover, it is easy to see that the base change map $\mathcal{B}_2^N(\Q,k) \rightarrow \mathcal{A}_{2,bc}^{(N)}(K,k)$ 
is {\em at most} two-to-one.

\begin{rem} Assume that $\ell=disc(K / \Q)$ is prime.  In this special case, the elliptic modular form associated to $\pi \in \mathcal{B}_2^N(\Q, k)$ lives in 
$S_k(\Gamma_0(N))^{new}$, $S_k(\Gamma_0(N\ell),\omega_K)^{new}$ and $S_K(\Gamma_0(N\ell^2),\epsilon)^{N\text{-}new}$ depending 
on whether $\pi_\ell$ has type $1,2$ or $3$ as above. The base change map is one-to-one on the subset of $\pi$'s with local type $1$ at $\ell$, and is two-to-one 
on the subset of $\pi's$ with local type $2$ or $3$ at $\ell$.
\end{rem}

To compute the dimensions of spaces of elliptic modular forms corresponding to $\pi$ with the local conditions as above, we use a formula (see \cite[Proposition 4.18.]{fgt}) 
that computes the multiplicity of any given representation of $\SL_2(\Z / N\Z)$ in the space $S_k(\Gamma(N))$ of elliptic modular forms of weight $k$ for the {\em principal} 
congruence subgroup of level $N$. 

\begin{prop} Let $N \ge 1$ and $k \ge 2$ be integers and $\sigma$ a representation of $G_N=\SL_2(\Z / N \Z)$ such that $\sigma(-I_2)$ is $(-1)^k$. 
Let $U_N \le G_N$ be the subgroup of all upper triangular unipotent elements and $S_3$ and $S_4$ be the images of elements of $\SL_2(\Z)$ of order 
$3$ and $4$, respectively. Then 
$$\Dim \Hom_{G_N}(\sigma, S_k(\Gamma(N))) = \dfrac{k-1}{12} \Dim \sigma - \dfrac{1}{2} \Dim \sigma^{U_N} + \varepsilon_k \tr \sigma(S_4) + \rho_k \tr \sigma(S_3) + 
\delta_{k,2} \Dim \sigma^{G_N}$$
where 
$$\varepsilon_k = \begin{cases}
\frac{(-1)^{k/2}}{4}, \ \ \ \textrm{if}\ \ \  k \equiv 0 \mod 2, \\
0, \ \ \ \ \ \ \  \ \ \ \ \textrm{otherwise}. \end{cases}$$
and 
$$\rho_k = \begin{cases}
0, \ \ \ \ \  \ \ \  \textrm{if}\ \ \  k \equiv 1 \mod 3, \\
-1/3, \ \ \ \textrm{if}\ \ \  k \equiv 2 \mod 3, \\
1/3, \ \ \ \ \ \textrm{if}\ \ \  k \equiv 0 \mod 3. \end{cases}$$
and $\delta_{k,2}$ is the Kronecker delta.
\end{prop} 

\begin{rem} Let us assume the situation of the previous remark. Then to compute the dimension of $S_k(\Gamma_0(N))^{new}$, we need to 
consider (the restriction to $SL_2(\Z/N\Z)$ of) the representation $\sigma$ of $\GL_2(\Z / N\Z)$ which is the tensor product of the representations 
$\sigma_p$ of $\GL_2(\F_p)$ where $\sigma_p$ is the (virtual) representation $[St] - [\mathbf{1}]$ where $St$ is the Steinberg representation. 
Using the character table of  $\GL_2(\F_p)$, we conclude that for {\bf even} $k \ge 2$,
$$\Dim S_k(\Gamma_0(N))^{new} = \dfrac{k-1}{12} \prod_{p | N} (p-1) + \varepsilon_k \prod_{p|N} \left ( \left(  \frac{-4}{p}\right ) -1 \right )
+ \rho_k \prod_{p|N} \left ( \left(  \frac{-3}{p}\right ) -1 \right ) + \delta_{k,2} (-1)^{div(N)}$$
where $div(N)$ is the number of prime divisors of $N$. 

Using Finis-Grunewald-Tirao \cite{fgt} and Bushnell-Henniart \cite{BushnellHenniart2006}, we also compute that for {\bf odd} $k \ge 3$,
\begin{center}
\begin{tabular}{rcl} 
$\dim S_k(\Gamma_0(N\ell),\omega_K)^{new}$&$=$&$\dfrac{k-1}{12} (\ell+1) \prod_{p | N} (p-1)$ \\
&&$+  \varepsilon_k \left ( \left(  \frac{-4}{\ell}\right ) +1 \right ) \cdot \prod_{p|N} \left ( \left(  \frac{-4}{p}\right ) -1 \right )$ \\
&&$+ \rho_k \left ( \left(  \frac{-3}{\ell}\right ) +1 \right ) \cdot \prod_{p|N} \left ( \left(  \frac{-3}{p}\right ) -1 \right )$ 
\end{tabular}
\end{center}

Finally, to compute the subspace of $S_K(\Gamma_0(N\ell^2),\epsilon)^{corr}$ of $S_K(\Gamma_0(N\ell^2),\epsilon)^{N\text{-}new}$ that have the local type $3$ at $\ell$, we proceed as above and find that 
for {\bf even} $k \ge 2$, 
\begin{center}
\begin{tabular}{rcl} 
$\dim S_k(\Gamma_0(N\ell^2),\epsilon)^{corr}$&$=$&$\dfrac{k-1}{12} (\ell-1) \prod_{p | N} (p-1)$ \\
&&$+  \varepsilon_k  \tr \sigma_{\ell}(S_4) \cdot \prod_{p|N} \left ( \left(  \frac{-4}{p}\right ) -1 \right )$ \\
&&$+ \rho_k \left ( \left(  \frac{-3}{\ell}\right ) -1 \right ) \cdot \prod_{p|N} \left ( \left(  \frac{-3}{p}\right ) -1 \right )$ 
\end{tabular}
\end{center}
where 
$$\tr \sigma_{\ell}(S_4) = \begin{cases} 2\cdot(-1)^{(\ell-3)/4}, \ \ \ \textrm{if}\ \ \  \ell \equiv 3 \mod 4, \\ 0, \ \ \ \ \ \ \ \ \ \ \ \ \ \ \ \  \ \ \ \ \textrm{otherwise}. \end{cases}$$

Let $BS^{bc}_k(\Gamma_0(N))^{new}$ denote the base-change subspace of space of {\em new} cuspidal {\em Bianchi modular forms} over $K$ of weight $k$ for $\Gamma_0(\langle N \rangle)$, that is, the subspace of $BS_k(\Gamma_0(N))^{new}$ that is accounted for by $\Pi$'s in $\mathcal{A}^N_{bc}(K,k)$. Then for {\bf even} $k \ge 2$, we have 
$$\Dim BS^{bc}_k(\Gamma_0(N))^{new}= \Dim S_k(\Gamma_0(N))^{new} + \dfrac{1}{2} \Dim S_k(\Gamma_0(N\ell^2),\epsilon)^{corr} $$
and for {\bf odd} $k \ge 3$, we have
$$\Dim BS^{bc}_k(\Gamma_0(N))^{new}= \dfrac{1}{2} \Dim S_k(\Gamma_0(N\ell),\omega_K)^{new}.$$

\end{rem}

\bibliography{BianchiBaseChange}
\bibliographystyle{amsplain}

\end{document}